\documentclass[a4paper,oneside,11pt]{article}%
\usepackage{amsmath}
\usepackage{amsfonts}
\usepackage{amssymb}
\usepackage{graphicx}%
\providecommand{\U}[1]{\protect\rule{.1in}{.1in}}
\newtheorem{theorem}{Theorem}[section]

\newtheorem{definition}[theorem]{Definition}
\newtheorem{assumption}[theorem]{Assumption}
\newtheorem{example}[theorem]{Example}

\newtheorem{lemma}[theorem]{Lemma}

\numberwithin{equation}{section}

\begin{document}

\title{Solutions for Functional Fully Coupled Forward-Backward Stochastic
Differential Equations}
\author{Shaolin Ji\thanks{Institute for Financial Studies and Institute of
Mathematics, Shandong University, Jinan, Shandong 250100, PR China
(Jsl@sdu.edu.cn, Fax: +86 0531 88564100).}
\and Shuzhen Yang\thanks{School of mathematics, Shandong University, Jinan,
Shandong 250100, PR China. (yangsz@mail.sdu.edu.cn). } }
\maketitle

\textbf{Abstract}: In this paper, we study a functional fully coupled
forward-backward stochastic differential equations (FBSDEs). Under a new type
of integral Lipschitz\ and {monotonicity} conditions, the existence and
uniqueness of solutions for functional fully coupled FBSDEs is proved. We also
investigate the relationship between the solution of functional fully coupled
FBSDE and the classical solution of the path-dependent {partial differential
equation} (P-PDE). When the solution of the P-PDE has some smooth and regular
properties, we solve the related functional fully coupled FBSDE and prove\ the
P-PDE\ has a\ unique solution.

\bigskip


\addcontentsline{toc}{section}{\hspace*{1.8em}Abstract}

\section{Introduction}

Linear backward stochastic differential equations (in short BSDEs) was
introduced by Bismut \cite{Bismut}. Pardoux and Peng \cite{Pardoux.E 3}
established the existence and uniqueness theorem for nonlinear BSDEs under a
standard Lipschitz condition. Since then, backward stochastic differential
equations and forward-backward stochastic differential equations (FBSDEs) have
been widely recognized that they provide useful tools in many fields,
especially mathematical finance and the stochastic control theory (see
\cite{Vitanic}, \cite{Dupire.B}, \cite{MY}, \cite{Yong1}, \cite{Yong3} and the
references therein).

A state dependent fully coupled FBSDE is formulated as:
\begin{equation}%
\begin{array}
[c]{cl}%
X(t)= & x+\int_{0}^{t}b(X(s),Y(s),Z(s))ds+\int_{0}^{t}\sigma
(X(s),Y(s),Z(s))dW(s),\\
Y(t)= & g(X(T))-\int_{t}^{T}h(X(s),Y(s),Z(s))ds-\int_{t}^{T}Z(s)dW(s),\quad
t\in\lbrack0,T].
\end{array}
\label{state1}%
\end{equation}
There have been three main methods to solve FBSDE (\ref{state1}), i.e., the
Method of Contraction Mapping (see \cite{Antonelli} and \cite{Pardoux.E}), the
Four Step Scheme (see \cite{Ma J.1}) and the Method of Continuation (see
\cite{Hu Y}, \cite{Peng} and \cite{Yong}). In \cite{Ma-J0}, Ma {\normalsize et
al.} studied the wellposedness of the FBSDEs in a general non-Markovian
framework. They find a unified scheme which combines all existing methodology
in the literature, and overcome some fundamental difficulties that have been
long-standing problems for non-Markovian FBSDEs.

In this paper, we study the following functional fully coupled FBSDE:%
\begin{equation}%
\begin{array}
[c]{cl}%
X(t)= & x+\int_{0}^{t}b(X_{s},Y(s),Z(s))ds+\int_{0}^{t}\sigma(X_{s}%
,Y(s),Z(s))dW(s),\\
Y(t)= & g(X(T))-\int_{t}^{T}h(X_{s},Y(s),Z(s))ds-\int_{t}^{T}Z(s)dW(s),\quad
t\in\lbrack0,T],
\end{array}
\label{state2}%
\end{equation}
where $X_{s}:=X(t)_{0\leq t\leq s}$.

As mentioned above, Hu and Peng \cite{Hu Y} initiated the continuation method
in which the key issue is a certain monotonicity condition. But unfortunately,
the Lipschitz and monotonicity\ conditions in \cite{Hu Y} and \cite{Peng} do
not work for equation (\ref{state2}). Here the main difficulty is that the
coefficients of (\ref{state2}) depend on the path of the solution $X(t)_{0\leq
t\leq T}$. In this paper, we revise the continuation method and propose a new
type of Lipschitz\ and monotonicity conditions. These new conditions involve
an integral term with respect to the path of $X(t)_{0\leq t\leq T}$. Thus, we
call them the integral Lipschitz\ and monotonicity conditions. The readers may
see Assumption \ref{ass-1} and \ref{ass-2} for more details. In particular, we
present two examples to illustrate that our assumptions are not restrictive.
Under the integral Lipschitz\ and monotonicity conditions, the continuation
method can go through and it leads to the existence and uniqueness of the
solution to equation (\ref{state2}).

It is well known that quasilinear parabolic partial differential equations are
related to Markovian FBSDEs (see \cite{Peng S 5}, \cite{Pardoux.E 2} and
\cite{Pardoux.E}), which generalizes the classical Feynman-Kac formula.
Recently a new framework of {functional It\^{o} calculus} was introduced by
Dupire \cite{Dupire.B} and later developed by Cont and {\normalsize Fourni
}\cite{Cont.R}, \cite{Cont-2}, \cite{Cont-3}. Inspired by Dupire's work, Peng
and Wang \cite{Peng S 3} obtained a so-called functional Feynman-Kac formula
for classical solutions of path-dependent {partial differential equation
(P-PDE)} in terms of non-Markovian BSDEs. Furthermore, under a special
condition, Peng {\cite{Peng S 2} proved that the viscosity solution of the
second order fully nonlinear P-PDE is unique.} Ekren, Touzi, and Zhang
({\cite{E}, \cite{E1}, \cite{E2})} gave another definition of the viscosity
solution of the fully nonlinear {P-PDE} and obtained the uniqueness result of
viscosity solutions.

In this paper, we explore the relationship between the solution of functional
fully coupled FBSDE (\ref{state2}) and the classical solution of the following
related P-PDE:
\begin{align*}
&  D_{t}u(\gamma_{t})+\mathcal{L}u(\gamma_{t})-h(\gamma_{t},u(\gamma
_{t}),v(\gamma_{t}))=0,\\
&  v(\gamma_{t})=D_{x}u(\gamma_{t})\tilde{\sigma}(\gamma_{t},u(\gamma
_{t}),v(\gamma_{t})),\\
&  u(\gamma_{T})=g(\gamma_{T}^{1},\gamma^{2}(T)),\quad\gamma_{T}^{1}\in
\Lambda^{d},\text{ }\gamma^{2}(T)\in\mathbb{R}^{n},
\end{align*}
where
\[
\mathcal{L}u=(\mathcal{L}u_{1},\cdots,\mathcal{L}u_{n}),\quad\mathcal{L}%
=\frac{1}{2}\sum_{i,j=1}^{n+d}(\tilde{\sigma}\tilde{\sigma}^{T})_{i,j}%
(\gamma_{t},u,v)D_{x_{i}x_{j}}+\sum_{i=1}^{n+d}\tilde{b}_{i}(\gamma
_{t},u,v)D_{x_{i}}.
\]
We prove that if the solution $u$ of the above P-PDE has some smooth and
regular properties, then we can solve the related equation (\ref{state2}) and
consequently,\ the P-PDE\ has a\ unique solution.

The paper is organized as follows. In section 2, we formulate the problem and
give the integral Lipschitz\ and monotonicity conditions. The existence and
uniqueness of the solution for (\ref{state2}) are proved in the first part of
section 3. In the second part of section 3, we show the relationship between
functional FBSDEs and P-PDEs.

\section{Formulation of the problem}

Let $\Omega=C([0,T];\mathbb{R}^{d})$ and $P$ the Wiener measure on
$(\Omega,\mathcal{B}(\Omega))$. We denote by $W=(W(t)_{t\in\lbrack0,T]})$ the
canonical Wiener process, with $W(t,\omega)=\omega(t)$, $t\in\lbrack0,T]$,
$\omega\in\Omega$. For any $t\in\lbrack0,T]$, we denote by $\mathcal{F}_{t}$
the $P$-completion of $\sigma(W(s),s\in\lbrack0,t])$.

For $n\in N$, set
\[
C_{t}^{n}=C([0,t];\mathbb{R}^{n})\text{ and }C^{n}=\bigcup_{t\in\lbrack
0,T]}C_{t}^{n}.
\]
Consider the following functional fully coupled FBSDE:
\begin{equation}
X(t)=x+\int_{0}^{t}b(X_{s},Y(s),Z(s))ds+\int_{0}^{t}\sigma(X_{s}%
,Y(s),Z(s))dW(s), \label{fbasd-1}%
\end{equation}%
\begin{equation}
Y(t)=g(X(T))-\int_{t}^{T}h(X_{s},Y(s),Z(s))ds-\int_{t}^{T}Z(s)dW(s),\quad
t\in\lbrack0,T], \label{fbase-2}%
\end{equation}
where the processes $X,Y,Z$ take values in $\mathbb{R}^{n},\mathbb{R}%
^{m},\mathbb{R}^{m\times d}$,\ $X_{s}=X(r)_{0\leq r\leq s}$\ and
\[%
\begin{array}
[c]{cl}%
b: & C^{n}\times\mathbb{R}^{m}\times\mathbb{R}^{m\times d}\times
\Omega\longrightarrow\mathbb{R}^{n};\\
\sigma: & C{^{n}}\times\mathbb{R}^{m}\times\mathbb{R}^{m\times d}\times
\Omega\longrightarrow\mathbb{R}^{n\times d};\\
h: & C{^{n}}\times\mathbb{R}^{m}\times\mathbb{R}^{m\times d}\times
\Omega\longrightarrow\mathbb{R}^{m};\\
g: & \mathbb{R}^{n}\times\Omega\longrightarrow\mathbb{R}^{m}.
\end{array}
\]
For $z\in\mathbb{R}^{m\times d}$, define $|z|=\{tr(zz^{T})\}^{1/2}$, where
\textquotedblleft$T$\textquotedblright\ means transpose. For $z^{1}%
\in\mathbb{R}^{m\times d}$, $z^{2}\in\mathbb{R}^{m\times d}$,%
\[
((z^{1},z^{2}))\triangleq tr(z^{1}(z^{2})^{T}).
\]
We use the notations
\[%
\begin{array}
[c]{rl}%
u^{1}= & (y^{1},z^{1})\in\mathbb{R}^{m}\times\mathbb{R}^{m\times d};\\
u^{2}= & (y^{2},z^{2})\in\mathbb{R}^{m}\times\mathbb{R}^{m\times d};\\
\lbrack u^{1},u^{2}]= & \langle y^{1},y^{2}\rangle+((z^{1},z^{2})).
\end{array}
\]
Given an $m\times n$ full-rank matrix $G,$ for $(x_{t},u)\in C{^{n}}%
\times\mathbb{R}^{m}\times\mathbb{R}^{m\times d}$, define%
\[
f(x_{t},u)=(G^{T}h(x_{t},u),Gb(x_{t},u),G\sigma(x_{t},u)),
\]
where $G\sigma=(G\sigma_{1},\cdots,G\sigma_{d}).$

\begin{definition}
We denote by $M^{2}(0,T;\mathbb{R}^{n})$ the set of all $\mathbb{R}^{n}%
$-valued $\mathcal{F}_{t}$-adapted processes $\vartheta(\cdot)$ such that
\[
E\int_{0}^{T}\mid\vartheta(s)\mid^{2}ds<+\infty.
\]

\end{definition}

\begin{definition}
A triple $(X,Y,Z):[0,T]\times\Omega\longrightarrow\mathbb{R}^{n}%
\times\mathbb{R}^{m}\times\mathbb{R}^{m\times d}$ is called an adapted
solution of the equations\ (\ref{fbasd-1}) and (\ref{fbase-2}), if $(X,Y,Z)\in
M^{2}(0,T;\mathbb{R}^{n}\times\mathbb{R}^{m}\times\mathbb{R}^{m\times d})$,
and it satisfies (\ref{fbasd-1}) and (\ref{fbase-2}) $P-a.s..$
\end{definition}

We rewrite (\ref{fbasd-1}) and (\ref{fbase-2}) in a differential form$:$
\begin{align*}
&  dX(t)=b(X_{t},Y(t),Z(t))dt+\sigma(X_{t},Y(t),Z(t))dW(t),\\
&  dY(t)=h(X_{t},Y(t),Z(t))dt+Z(t)dW(t),\\
&  X(0)=x,\ \ \ Y(T)=g(X(T)).
\end{align*}

Now we give the assumptions:

\begin{assumption}
\label{ass-1} For each $(x_{T},u_{T})\in C^{n}\times C^{m}\times C^{m\times
d}$, $f(x_{\cdot},u(\cdot))\in M^{2}(0,T;\mathbb{R}^{n}\times\mathbb{R}%
^{m}\times\mathbb{R}^{m\times d}).$ There exists a constant $c_{1}>0$, such
that%
\[%
\begin{array}
[c]{l}%
\int_{0}^{T}\left\vert f(x_{t}^{1},u^{1}(t))-f(x_{t}^{2},u^{1}(t))\right\vert
^{2}dt\leq c_{1}\int_{0}^{T}\left\vert x^{1}(t)-x^{2}(t)\right\vert ^{2}dt,\\
\left\vert f(x_{s}^{1},u^{1}(s))-f(x_{s}^{1},u^{2}(s))\right\vert \leq
c_{1}\left\vert u^{1}(s)-u^{2}(s)\right\vert ,\text{ \ \ }0\leq s\leq T,
\end{array}
\]
$\forall(x_{T}^{1},u_{T}^{1}),(x_{T}^{2},u_{T}^{2})\in C^{n}\times C^{m}\times
C^{m\times d},$ and
\[
\left\vert g(x^{1})-g(x^{2})\right\vert \leq c_{1}\left\vert x^{1}%
-x^{2}\right\vert ,\forall(x^{1},x^{2})\in\mathbb{R}^{n}{\times}\mathbb{R}%
^{n},\text{ \ \ }P-a.s..
\]

\end{assumption}

\begin{assumption}
\label{ass-2} There exists constants $\beta_{1},\beta_{2},\mu_{1}$ such that
\[%
\begin{array}
[c]{cl}
& \int_{0}^{T}[f(x_{t}^{1},u^{1}(t))-f(x_{t}^{2},u^{2}(t)),(x^{1}%
(t)-x^{2}(t),u^{1}(t)-u^{2}(t))]dt\\
\geq & \int_{0}^{T}[\beta_{1}\left\vert Gx^{1}(t)-Gx^{2}(t)\right\vert
^{2}+\beta_{2}(\left\vert G^{T}y^{1}(t)-G^{T}y^{2}(t)\right\vert
^{2}+\left\vert G^{T}z^{1}(t)-G^{T}z^{2}(t)\right\vert ^{2})]dt,\text{
\ \ }P-a.s.,
\end{array}
\]
\newline$\forall(x_{T}^{1},u_{T}^{1}),(x_{T}^{2},u_{T}^{2})\in C^{n}\times
C^{m}\times C^{m\times d},$ and
\[
\langle g(x^{1})-g(x^{2}),G(x^{1}-x^{2})\rangle\leq-\mu_{1}\left\vert
Gx^{1}-Gx^{2}\right\vert ^{2},\quad\ \ P-a.s.,
\]
$\forall(x^{1},x^{2})\in\mathbb{R}^{n}{\times}\mathbb{R}^{n},$ where
$\beta_{1},\beta_{2}$ and $\mu_{1}$ are given nonnegative constants with
$\beta_{1}+\beta_{2}>0,\mu_{1}+\beta_{2}>0.$ Moreover we have $\beta_{1}%
>0,\mu_{1}>0$ (resp., $\beta_{2}>0,\mu_{1}>0$) when $m>n$ (resp., $n>m$).
\end{assumption}

The coefficients of the following two examples satisfy Assumptions \ref{ass-1}
and \ref{ass-2}.

\begin{example}
\label{example}Set $n=m=1,$ $T=1,$ and the functional FBSDE is%
\begin{align}
&  dX(t)=(\int_{0}^{t}X(s)ds+2Y(t))dt+(\int_{0}^{t}%
X(s)ds+2Z(t))dW(t),\label{ex2}\\
&  dY(t)=(\int_{0}^{t}X(s)ds+3X(t))dt+Z(t)dW(t),\nonumber\\
&  X(0)=x,\ \ \ Y(1)=-X(1).\nonumber
\end{align}

We first show that the coefficients of (\ref{ex2}) satisfy Assumption
\ref{ass-1}. For any $(x_{1}^{1},y_{1}^{1},z_{1}^{1}),(x_{1}^{2},y_{1}%
^{2},z_{1}^{2})\in C\times C\times C{,}$ it is easy to check that
$g$\ satisfies the Lipschitz condition. For $f=(h,b,\sigma)$, we have%
\[%
\begin{array}
[c]{rl}
& \int_{0}^{1}\left\vert f(x_{t}^{1},y^{1}(t),z^{1}(t))-f(x_{t}^{2}%
,y^{1}(t),z^{1}(t))\right\vert ^{2}dt\\
\leq & \int_{0}^{1}[4t^{2}\int_{0}^{t}\left\vert x^{1}(s)-x^{2}(s)\right\vert
^{2}ds+18\left\vert x^{1}(t)-x^{2}(t)\right\vert ^{2}]dt\\
\leq & 22\int_{0}^{1}\left\vert x^{1}(t)-x^{2}(t)\right\vert ^{2}dt.
\end{array}
\]

and%
\[%
\begin{array}
[c]{rl}
& \left\vert f(x_{t}^{1},y^{1}(t),z^{1}(t))-f(x_{t}^{1},y^{2}(t),z^{2}%
(t))\right\vert \\
\leq & 2\left\vert y^{1}(t)-y^{2}(t)\right\vert +2\left\vert z^{1}%
(t)-z^{2}(t)\right\vert .
\end{array}
\]

Then we check Assumption \ref{ass-2}. For any $(x_{1}^{1},y_{1}^{1},z_{1}%
^{1}),(x_{1}^{2},y_{1}^{2},z_{1}^{2})\in C\times C\times C{,}$ $g$ satisfies
the monotonicity condition:%
\[
(g(x^{1}(1))-g(x^{2}(1)))\cdot(x^{1}(1)-x^{2}(1))\leq-\left\vert
x^{1}(1)-x^{2}(1)\right\vert ^{2}.
\]
For $f=(h,b,\sigma)$, we have%
\[%
\begin{array}
[c]{cl}
& \int_{0}^{1}[f(x_{t}^{1},u^{1}(t))-f(x_{t}^{2},u^{2}(t)),(x^{1}%
(t)-x^{2}(t),u^{1}(t)-u^{2}(t))]dt\\
= & \int_{0}^{1}[\int_{0}^{t}(x^{1}(s)-x^{2}(s))ds(x^{1}(t)-x^{2}%
(t))+3\left\vert x^{1}(t)-x^{2}(t)\right\vert ^{2}\\
& +(\int_{0}^{t}x^{1}(s)-x^{2}(s)ds)\cdot(y^{1}(t)-y^{2}(t))+2\left\vert
y^{1}(t)-y^{2}(t)\right\vert ^{2}\\
& +(\int_{0}^{t}x^{1}(s)-x^{2}(s)ds)\cdot(z^{1}(t)-z^{2}(t))+2\left\vert
z^{1}(t)-z^{2}(t)\right\vert ^{2}]dt\\
\geq & \int_{0}^{1}[\left\vert x^{1}(t)-x^{2}(t)\right\vert ^{2}+\left\vert
y^{1}(t)-y^{2}(t)\right\vert ^{2}+\left\vert z^{1}(t)-z^{2}(t)\right\vert
^{2}]dt.
\end{array}
\]

\end{example}

\begin{example}
We consider a nonlinear functional fully coupled FBSDE. Set $n=m=1,$ $T=1$.
Suppose that $\hat{f}=(\hat{h},\hat{b},\hat{\sigma})$ satisfies the next
Lipschitz and monotonicity conditions:

there exists a constant $c_{1}>0$, such that%
\[
\left\vert \hat{f}(x^{1},y^{1},z^{1})-\hat{f}(x^{2},y^{2},z^{2})\right\vert
\leq c_{1}(\left\vert x^{1}-x^{2}\right\vert +\left\vert y^{1}-y^{2}%
\right\vert +\left\vert z^{1}-z^{2}\right\vert ),\ \ \ P-a.s.,
\]

\noindent and
\[%
\begin{array}
[c]{cl}
& [\hat{f}(x^{1},y^{1},z^{1})-\hat{f}(x^{2},y^{2},z^{2}),(x^{1}-x^{2}%
,y^{1}-y^{2},z^{1}-z^{2})]\\
\geq & c_{1}(\left\vert x^{1}-x^{2}\right\vert ^{2}+\left\vert y^{1}%
-y^{2}\right\vert ^{2}+\left\vert z^{1}-z^{2}\right\vert ^{2}),\ \ \ P-a.s.,
\end{array}
\]
\newline$\forall(x^{1},y^{1},z^{1}),(x^{2},y^{2},z^{2})\in\mathbb{R}%
\times\mathbb{R}\times\mathbb{R}$. Then we set
\[%
\begin{array}
[c]{rl}%
h(x_{t},y,z)= & \hat{h}(\int_{0}^{t}x(s)ds,y,z);\\
b(x_{t},y,z)= & \hat{b}(\int_{0}^{t}x(s)ds,y,z);\\
\sigma(x_{t},y,z)= & \hat{\sigma}(\int_{0}^{t}x(s)ds,y,z);\\
g(x(1))= & -x(1),
\end{array}
\]
$\forall(x_{1},y,z)\in C\times\mathbb{R}\times\mathbb{R}$. The corresponding
functional FBSDE is%
\begin{align*}
&  dX(t)=\hat{b}(\int_{0}^{t}X(s)ds,Y(t),Z(t))dt+\hat{\sigma}(\int_{0}%
^{t}X(s)ds,Y(t),Z(t))dW(t),\\
&  dY(t)=\hat{h}(\int_{0}^{t}X(s)ds,Y(t),Z(t))dt+Z(t)dW(t),\\
&  X(0)=x,\ \ \ Y(1)=-X(1).
\end{align*}

Applying similar analysis as in Example \ref{example}, we can check that
$f=(h,b,\sigma)$ and$\ g$ satisfy Assumptions \ref{ass-1} and \ref{ass-2}.
\end{example}

\section{Functional Fully Coupled FBSDEs}

\subsection{Existence and uniqueness for functional\ fully coupled FBSDEs}

In this section, we combine the continuation method in \cite{Peng} and the
integral Lipschitz\ and monotonicity conditions to prove our main result.

\begin{theorem}
\label{theorem-1} Let Assumptions \ref{ass-1} and \ref{ass-2} hold. Then there
exist a unique adapted solution $(X,Y,Z)$ for equations (\ref{fbasd-1})and
(\ref{fbase-2}).
\end{theorem}

Firstly, we prove the uniqueness.

\textbf{Proof}. Let $U^{1}(\cdot)\triangleq(Y^{1}(\cdot),Z^{1}(\cdot))$,
$U^{2}(\cdot)\triangleq(Y^{2}(\cdot),Z^{2}(\cdot))$ be two adapted solutions
of (\ref{fbasd-1})and (\ref{fbase-2}). We set
\begin{align*}
(\hat{X},\hat{Y},\hat{Z})  &  =(X^{1}-X^{2},Y^{1}-Y^{2},Z^{1}-Z^{2});\\
\hat{b}(t)  &  =b(X_{t}^{1},U^{1}(t))-b(X_{t}^{2},U^{2}(t));\\
\hat{\sigma}(t)  &  =\sigma(X_{t}^{1},U^{1}(t))-\sigma(X_{t}^{2},U^{2}(t));\\
\hat{h}(t)  &  =h(X_{t}^{1},U^{1}(t))-h(X_{t}^{2},U^{2}(t)).
\end{align*}
From Assumption \ref{ass-1}, it follows that $\hat{X}(\cdot)$ and $\hat
{Y}(\cdot)$ are continuous, and
\[
E[\int_{0}^{T}(\left\vert \hat{X}(t)\right\vert ^{2}+\left\vert \hat
{Y}(t)\right\vert ^{2})dt]<+\infty.
\]
Applying the It\^{o} formula to $\langle\hat{Y}(t),G\hat{X}(t)\rangle$,%
\[%
\begin{array}
[c]{cl}%
d\langle\hat{Y}(t),G\hat{X}(t)\rangle= & [f(X_{t}^{1},U^{1}(t))-f(X_{t}%
^{2},U^{2}(t)),(X^{1}(t)-X^{2}(t),U^{1}(t)-U^{2}(t))]dt\\
& +(\hat{X}^{T}(t)G^{T}\hat{Z}^{T}(t)+\hat{Y}^{T}(t)G\hat{\sigma}(t))dW(t).
\end{array}
\]
It yields that%
\[%
\begin{array}
[c]{cl}
& E\langle g(X^{1}(T))-g(X^{2}(T)),GX^{1}(T)-GX^{2}(T)\rangle\\
= & E\int_{0}^{T}[f(X_{t}^{1},U^{1}(t))-f(X_{t}^{2},U^{2}(t)),(X^{1}%
(t)-X^{2}(t),U^{1}(t)-U^{2}(t))]dt.
\end{array}
\]
\ \ \ By Assumptions \ref{ass-1} and \ref{ass-2}, we obtain%
\[%
\begin{array}
[c]{cl}
& -\mu_{1}E\mid GX^{1}(T)-GX^{2}(T)\mid^{2}\\
\geq & E\langle g(X^{1}(T))-g(X^{2}(T)),GX^{1}(T)-GX^{2}(T)\rangle\\
\geq & E\int_{0}^{T}(\beta_{1}\mid Gx^{1}(t)-Gx^{2}(t)\mid^{2}+\beta_{2}(\mid
G^{T}Y^{1}(t)-G^{T}Y^{2}(t)\mid+\mid G^{T}Z^{1}(t)-G^{T}Z^{2}(t)\mid))dt.
\end{array}
\]

In the case $m>n$, we have $\beta_{1}>0$ and $\mu_{1}>0$. Then we have
$X^{1}=X^{2}.$ From the uniqueness of BSDE, it follows that $U^{1}=U^{2}.$ In
the case $m<n$, we have $\beta_{2}>0$ and $\mu_{1}>0$. It yields that
$U^{1}=U^{2}.$ Thus, by the uniqueness of SDE, it yields that $X^{1}=X^{2}$.
$\ \ \Box$

The following\ lemma is quoted from \cite{Peng}.

\begin{lemma}
Suppose that $(b_{0}(\cdot),\sigma_{0}(\cdot),h_{0}(\cdot))\in M^{2}%
(0,T;\mathbb{R}^{n}\times\mathbb{R}^{n\times d}\times\mathbb{R}^{m}),$
$g_{0}\in L^{2}(\Omega,\mathcal{F}_{T};\mathbb{R}^{m})$. Then the following
linear forward-backward stochastic differential equation
\begin{equation}
X(t)=x+\int_{0}^{t}(-\beta_{2}G^{T}Y(s)+b_{0}(s))ds+\int_{0}^{t}(-\beta
_{2}G^{T}Z(s)+\sigma_{0}(s))dW(s), \label{lemma-eqs-1}%
\end{equation}%
\begin{equation}
Y(t)=(\lambda GX(T)+g_{0})-\int_{t}^{T}(-\beta_{1}GX(s)+h_{0}(s))ds-\int%
_{t}^{T}Z(s)dW(s) \label{lemma-eqs-2}%
\end{equation}
have a unique adapted solution $(X,Y,Z),$ where $\lambda$ is a nonnegative constant.
\end{lemma}

For any given $\alpha\in\mathbb{R}$, we define
\begin{align*}
b^{\alpha}(x_{t},y,z)  &  =\alpha b(x_{t},y,z)+(1-\alpha)\beta_{2}(-G^{T}y);\\
\sigma^{\alpha}(x_{t},y,z)  &  =\alpha\sigma(x_{t},y,z)+(1-\alpha)\beta
_{2}(-G^{T}z);\\
h^{\alpha}(x_{t},y,z)  &  =\alpha h(x_{t},y,z)+(1-\alpha)\beta_{1}(-Gx);\\
g^{\alpha}(x)  &  =\alpha g(x)+(1-\alpha)\beta_{1}(Gx),
\end{align*}
and consider the following equations:
\begin{equation}
X(t)=x+\int_{0}^{t}(b^{\alpha}(X_{s},U(s))+b_{0}(s))ds+\int_{0}^{t}%
(\sigma^{\alpha}(X_{s},U(s))+\sigma_{0}(s))dW(s), \label{eqs-3}%
\end{equation}%
\begin{equation}
Y(t)=(g^{\alpha}(X(T))+g_{0})-\int_{t}^{T}(h^{\alpha}(X_{s},U(s))+h_{0}%
(s))ds-\int_{t}^{T}Z(s)dW(s), \label{eqs-4}%
\end{equation}
where $U=(Y,Z)$.

\begin{lemma}
\label{lemma-3.2} We assume that, a priori, for a given $\alpha_{0}\in
\lbrack0,1)$ and for any $(b_{0}(\cdot),\sigma_{0}(\cdot),h_{0}(\cdot))\in
M^{2}(0,T;\mathbb{R}^{n}\times\mathbb{R}^{n\times d}\times\mathbb{R}%
^{m}),\quad g_{0}\in L^{2}(\Omega,\mathcal{F}_{T};\mathbb{R}^{m})$, there
exists a solution of (\ref{eqs-3}) and (\ref{eqs-4}). Then there exists a
$\delta_{0}\in(0,1)$ which depends only on $\mu_{1},\beta_{1},\beta_{2}$ and
$T$, such that for all $\alpha\in\lbrack\alpha_{0},\alpha_{0}+\delta_{0}]$,
and for any $(b_{0}(\cdot),\sigma_{0}(\cdot),h_{0}(\cdot))\in M^{2}%
(0,T;\mathbb{R}^{n}\times\mathbb{R}^{n\times d}\times\mathbb{R}^{m}),\quad
g_{0}\in L^{2}(\Omega,\mathcal{F}_{T};\mathbb{R}^{m})$,\ (\ref{eqs-3}) and
(\ref{eqs-4}) have an adapted solution.
\end{lemma}

\textbf{Proof}. Note that%
\begin{align*}
b^{\alpha_{0}+\delta}(x_{t},y,z)  &  =b^{\alpha_{0}}(x_{t},y,z)+\delta
(\beta_{2}G^{T}y+b(x_{t},y,z));\\
\sigma^{\alpha_{0}+\delta}(x_{t},y,z)  &  =\sigma^{\alpha_{0}}(x_{t}%
,y,z)+\delta(\beta_{2}G^{T}z+\sigma(x_{t},y,z));\\
h^{\alpha_{0}+\delta}(x_{t},y,z)  &  =h^{\alpha_{0}}(x_{t},y,z)+\delta
(\beta_{1}Gx+h(x_{t},y,z));\\
g^{\alpha_{0}+\delta}(x)  &  =g^{\alpha_{0}}(x)+\delta(-\beta_{1}Gx+g(x)).
\end{align*}
\ Let $(X^{0},U^{0})=(X^{0},Y^{0},Z^{0})=0$. We solve the following equations
iteratively:%
\begin{equation}%
\begin{array}
[c]{cl}%
X^{i+1}(t)= & x+\int_{0}^{t}[b^{\alpha_{0}}(X_{s}^{i+1},U^{i+1}(s))+\delta
(\beta_{2}G^{T}Y^{i}(s)+b(X_{s}^{i},U^{i}(s)))+b_{0}(s)]ds\\
& +\int_{0}^{t}[\sigma^{\alpha_{0}}(X_{s}^{i+1},U^{i+1}(s))+\delta(\beta
_{2}G^{T}Z^{i}(s)+\sigma(X_{s}^{i},U^{i}(s)))+\sigma_{0}(s)]dW(s),
\end{array}
\label{3.5}%
\end{equation}%
\begin{equation}%
\begin{array}
[c]{cl}%
Y^{i+1}(t)= & [g^{\alpha_{0}}(X^{i+1}(T))+\delta(-\beta_{1}GX^{i}%
(T)+g(X^{i}(T)))+g_{0}(s)]-\int_{t}^{T}Z^{i+1}(s)dW(s)\\
& -\int_{t}^{T}[h^{\alpha_{0}}(X_{s}^{i+1},U^{i+1}(s))+\delta(\beta_{1}%
GX^{i}(s)+h(X_{s}^{i},U^{i}(s)))+h_{0}(s)]ds,
\end{array}
\label{3.6}%
\end{equation}
where $U^{i}=(Y^{i},Z^{i})$.

Set
\[
(\hat{X}^{i+1},\hat{U}^{i+1})=(\hat{X}^{i+1},\hat{Y}^{i+1},\hat{Z}^{i+1}).
\]
Applying It\^{o} formula to $\langle\hat{Y}^{i+1}(t),G\hat{X}^{i+1}(t)\rangle
$,%
\[%
\begin{array}
[c]{cl}
& E\langle g^{\alpha_{0}}(X^{i+1}(T))-g^{\alpha_{0}}(X^{i}(T)),G\hat{X}%
^{i+1}(T)\rangle\\
= & \delta E\langle\beta_{1}G\hat{X}^{i}(T)-(g(X^{i}(T))-g(X^{i-1}%
(T))),G\hat{X}^{i+1}(T)\rangle\\
& +E\int_{0}^{T}[f^{\alpha_{0}}(X_{t}^{i+1},U^{i+1}(t))-f^{\alpha_{0}}%
(X_{t}^{i},U^{i}(t)),(\hat{X}^{i+1}(t),\hat{U}^{i+1}(t))]dt\\
& +\delta(E\int_{0}^{T}[(\beta_{1}G^{T}G\hat{X}^{i}(t),\beta_{2}GG^{T}\hat
{U}^{i}(t))+f(X_{t}^{i},U^{i}(t))-f(X_{t}^{i-1},U^{i-1}(t)),(\hat{X}%
^{i+1}(t),\hat{U}^{i+1}(t))]).
\end{array}
\]
From Assumptions \ref{lemma-eqs-1} and \ref{lemma-eqs-2}, there exists
$K_{1}>0$ such that
\begin{equation}%
\begin{array}
[c]{cl}
& (\alpha_{0}\mu_{1}+(1-\alpha_{0}))E\mid G\hat{X}^{i+1}(T)\mid^{2}+\beta
_{1}E\int_{0}^{T}\mathbb{\mid}G\hat{X}^{i+1}(t)\mathbb{\mid}^{2}dt+\beta
_{2}E\int_{0}^{T}\mid G^{T}\hat{U}^{i+1}(t)\mid^{2}dt\\
\leq & \delta K_{1}(E\mid\hat{X}^{i}(T)\mid^{2}+E\int_{0}^{T}\mathbb{\mid}%
\hat{X}^{i}(t)\mathbb{\mid}^{2}dt+E\int_{0}^{T}\mid\hat{U}^{i}(t)\mid^{2}dt\\
& +E\mid\hat{X}^{i+1}(T)\mid^{2}+E\int_{0}^{T}\mathbb{\mid}\hat{X}%
^{i+1}(t))\mathbb{\mid}^{2}dt+E\int_{0}^{T}\mid\hat{U}^{i+1}(t)\mid^{2}dt).
\end{array}
\label{FB-1}%
\end{equation}
We also have $\forall i\geq1,$%
\[%
\begin{array}
[c]{cl}
& \hat{X}^{i+1}(s)\\
= & \int_{0}^{s}[b^{\alpha_{0}}(X_{t}^{i+1},U^{i+1}(t))-b^{\alpha_{0}}%
(X_{t}^{i},U^{i}(t))+\delta(\hat{Y}^{i}(t)+b(X_{t}^{i},U^{i}(t))-b(X_{t}%
^{i-1},U^{i-1}(t)))]dt\\
& +\int_{0}^{s}[\sigma^{\alpha_{0}}(X_{t}^{i+1},U^{i+1}(t))-\sigma^{\alpha
_{0}}(X_{t}^{i},U^{i}(t))+\delta(\hat{Z}^{i}(t)+\sigma(X_{t}^{i}%
,U^{i}(t))-\sigma(X_{t}^{i-1},U^{i-1}(t)))]dW(t).
\end{array}
\]
By the usual technique of\ estimation, we can derive that
\begin{equation}%
\begin{array}
[c]{cl}
& \sup_{0\leq s\leq T}E\mid\hat{X}^{i+1}(s)\mid^{2}\\
\leq & K_{1}\delta(E\int_{0}^{T}\mid\hat{U}^{i}(t)\mid^{2}dt+E\int_{0}%
^{T}\mathbb{\mid}\hat{X}^{i}(t)\mathbb{\mid}^{2}dt)+K_{1}E\int_{0}^{T}\mid
\hat{U}^{i+1}(t)\mid^{2}dt,\quad\ \forall i\geq1.
\end{array}
\label{FB-2}%
\end{equation}

and%

\begin{equation}%
\begin{array}
[c]{cl}
& E\int_{0}^{T}\mathbb{\mid}\hat{X}^{i+1}(t)\mathbb{\mid}^{2}dt\\
\leq & K_{1}T\delta(E\int_{0}^{T}\mid\hat{U}^{i}(t)\mid^{2}dt+E\int_{0}%
^{T}\mathbb{\mid}\hat{X}^{i}(t)\mathbb{\mid}^{2}dt)+K_{1}TE\int_{0}^{T}%
\mid\hat{U}^{i+1}(t)\mid^{2}dt,\quad\ \forall i\geq1.
\end{array}
\label{FB-3}%
\end{equation}
For equation (\ref{3.6}), by the standard estimation of BSDE, we have%

\begin{equation}%
\begin{array}
[c]{cl}
& E\int_{0}^{T}\mid\hat{U}^{i+1}(t)\mid^{2}dt\\
\leq & K_{1}\delta(E\int_{0}^{T}\mid\hat{U}^{i}(t)\mid^{2}dt+E\int_{0}%
^{T}\mathbb{\mid}\hat{X}^{i}(t)\mathbb{\mid}^{2}dt+K_{1}E\mid\hat{X}%
^{i}(T)\mid^{2})\\
& +K_{1}E\int_{0}^{T}\mathbb{\mid}\hat{X}^{i+1}(t)\mathbb{\mid}^{2}%
dt+K_{1}E\mid\hat{X}^{i+1}(T)\mid^{2},
\end{array}
\label{FB-4}%
\end{equation}
where the constant $K_{1}$ depends on the Lipschitz constants $c_{1},\beta
_{1},\beta_{2}$ and $T.$

Due to the above estimations, we conclude that
\[%
\begin{array}
[c]{cl}
& E\int_{0}^{T}\mid\hat{U}^{i+1}(t)\mid^{2}dt+E\mid\hat{X}^{i+1}(T)\mid
^{2}+E\int_{0}^{T}\mathbb{\mid}\hat{X}^{i+1}(t)\mathbb{\mid}^{2}dt\\
\leq & K\delta(E\int_{0}^{T}\mid\hat{U}^{i}(t)\mid^{2}dt+E\mid\hat{X}%
^{i}(T)\mid^{2}+E\int_{0}^{T}\mathbb{\mid}\hat{X}^{i}(t)\mathbb{\mid}^{2}dt),
\end{array}
\]
where the constant $K$\ depends on $\beta_{1},$ $\beta_{2},$ $\mu,K_{1}$ and
$T.$ Taking $\delta_{0}=\frac{1}{2K},$ for each $0<\delta\leq\delta_{0},$%
\[%
\begin{array}
[c]{cl}
& E\int_{0}^{T}\mid\hat{U}^{i+1}(t)\mid^{2}dt+E\mid\hat{X}^{i+1}(T)\mid
^{2}+E\int_{0}^{T}\mathbb{\mid}\hat{X}^{i+1}(t)\mathbb{\mid}^{2}dt\\
\leq & \frac{1}{2}(E\int_{0}^{T}\mid\hat{U}^{i}(t)\mid^{2}dt+E\mid\hat{X}%
^{i}(T)\mid^{2}+E\int_{0}^{T}\mathbb{\mid}\hat{X}^{i}(t)\mathbb{\mid}^{2}dt).
\end{array}
\]

Thus, the limit of $(U^{i},X^{i})$ exists. We denote its limit by
$(X,U)=(X,Y,Z)$. Passing to the limit in equations (\ref{3.5}) and
(\ref{3.6}), when $0<\delta\leq\delta_{0}$, $(X,U)=(X,Y,Z)$ solves equations
(\ref{eqs-3}) and (\ref{eqs-4}) for $\alpha=\alpha_{0}+\delta$. This completes
the proof. $\ \ \Box$

In the following, we give the proof of the existence.

\textbf{Proof}. By Lemma 3.1, when $\alpha=0$, for any $(b_{0}(\cdot
),\sigma_{0}(\cdot),h_{0}(\cdot))\in M^{2}(0,T;\mathbb{R}^{n}\times
\mathbb{R}^{n\times d}\times\mathbb{R}^{m}),\quad g_{0}\in L^{2}%
(\Omega,\mathcal{F}_{T};\mathbb{R}^{m})$ equations (\ref{eqs-3}) and
(\ref{eqs-4}) have an adapted solution. Then from Lemma \ref{lemma-3.2}, there
exists a constant $\delta_{0}$ which depends only on $\mu_{1},\beta_{1}%
,\beta_{2}$ and $T$, such that for all $\alpha\in\lbrack\alpha_{0},\alpha
_{0}+\delta_{0}]$,\ (\ref{eqs-3}) and (\ref{eqs-4}) have an adapted solution.
Thus, we can solve equations (\ref{eqs-3}) and (\ref{eqs-4}) successively for
the case $\alpha\in\lbrack0,\delta_{0}],[\delta_{0},2\delta_{0}],\ldots$ which
leads to that, when $\alpha=1$, there exists an adapted solution of equations
(\ref{fbasd-1}) and (\ref{fbase-2}). $\ \ \Box$

\subsection{Relationship between functional fully coupled FBSDEs
(\ref{state2}) and related P-PDEs}

In this section, we first review some basic notions and results of {functional
It\^{o} calculus. Then we investigate }the relationship between the solution
of functional fully coupled FBSDE and the classical solution of the related P-PDE.

Let $T>0$ be fixed. For each $t\in\lbrack0,T]$, we denote by $\Lambda_{t}^{n}$
the set of c\`{a}dl\`{a}g $\mathbb{R}^{n}$-valued functions on $[0,t]$.

For each $\gamma\in\Lambda_{T}^{n}$ the value of $\gamma$ at time $s\in
\lbrack0,T]$ is denoted by $\gamma(s)$. Thus $\gamma=\gamma(s)_{0\leq s\leq
T}$ is a c\`{a}dl\`{a}g process on $[0,T]$ and its value at time $s$ is
$\gamma(s)$. The path of $\gamma$ up to time $t$ is denoted by $\gamma_{t}$,
i.e., $\gamma_{t}=\gamma(s)_{0\leq s\leq t}\in\Lambda_{t}^{n}$. We denote
$\Lambda^{n}=\bigcup_{t\in\lbrack0,T]}\Lambda_{t}^{n}$. For each $\gamma
_{t}\in\Lambda^{n}$ and $x\in\mathbb{R}^{n}$ we denote by $\gamma_{t}(s)$ the
value of $\gamma_{t}$ at $s\in\lbrack0,t]$ and $\gamma_{t}^{x}:=(\gamma
_{t}(s)_{0\leq s<t},\gamma_{t}(t)+x)$ which is also an element in $\Lambda
_{t}^{n}$.

We now introduce a distance on $\Lambda^{n}$. Let $\langle\cdot,\cdot\rangle$
and $|\cdot|$ denote the inner product and norm in $\mathbb{R}^{n}$. For each
$0\leq t,\bar{t}\leq T$ and $\gamma_{t},\bar{\gamma}_{\bar{t}}\in\Lambda^{n}$,
we denote%
\[%
\begin{array}
[c]{l}%
\Vert\gamma_{t}\Vert:=\sup\limits_{s\in\lbrack0,t]}|\gamma_{t}(s)|,\\
\Vert\gamma_{t}-\bar{\gamma}_{\bar{t}}\Vert:=\sup\limits_{s\in\lbrack
0,t\vee\bar{t}]}|\gamma_{t}(s\wedge t)-\bar{\gamma}_{\bar{t}}(s\wedge\bar
{t})|,\\
d_{\infty}(\gamma_{t},\bar{\gamma}_{\bar{t}}):=\sup\limits_{s\in\lbrack
0,t\vee\bar{t}]}|\gamma_{t}(s\wedge t)-\bar{\gamma}_{\bar{t}}(s\wedge\bar
{t})|+|t-\bar{t}|.
\end{array}
\]

It is obvious that $\Lambda_{t}^{n}$ is a Banach space with respect to
$\Vert\cdot\Vert$. Since $\Lambda^{n}$ is not a linear space, $d_{\infty}$ is
not a norm.

\begin{definition}
A function $u:\Lambda^{n}\mapsto\mathbb{R}$ is said to be $\Lambda^{n}%
$--continuous at $\gamma_{t}\in\Lambda^{n}$, if for any $\varepsilon>0$ there
exists $\delta>0$ such that for each $\bar{\gamma}_{\bar{t}}\in\Lambda^{n}$
with $d_{\infty}(\gamma_{t},\bar{\gamma}_{\bar{t}})<\delta$, we have
$|u(\gamma_{t})-u(\bar{\gamma}_{\bar{t}})|<\varepsilon$.
\end{definition}

$u$ is said to be $\Lambda^{n}$--continuous if it is $\Lambda^{n}$--continuous
at each $\gamma_{t}\in\Lambda^{n}$.

\begin{definition}
Let $v:\Lambda^{n}\mapsto\mathbb{R}$ and $\gamma_{t}\in\Lambda^{n}$ be given.
If there exists $p\in\mathbb{R}^{n}$, such that
\[
v(\gamma_{t}^{x})=v(\gamma_{t})+\langle p,x\rangle+o(|x|)\ \text{as}%
\ x\rightarrow0,\ x\in\mathbb{R}^{n}.\ \
\]

Then we say that $v$ is (vertically) differentiable at $\gamma_{t}$ and denote
the gradient of $D_{x}v(\gamma_{t})=p$. $v$ is said to be vertically
differentiable in $\Lambda^{n}$ if $D_{x}v(\gamma_{t})$ exists for each
$\gamma_{t}\in\Lambda^{n}$. We can similarly define the Hessian $D_{xx}%
^{2}v(\gamma_{t})$. It is an $\mathbb{S}(n)$-valued function defined on
$\Lambda^{n}$, where $\mathbb{S}(n)$ is the space of all $n\times n$ symmetric matrices.
\end{definition}

For each $\gamma_{t}\in\Lambda^{n}$ we denote
\[
\gamma_{t,s}(r)=\gamma_{t}(r)\mathbf{1}_{[0,t)}(r)+\gamma_{t}(t)\mathbf{1}%
_{[t,s]}(r),\ \ r\in\lbrack0,s].
\]
It is clear that $\gamma_{t,s}\in\Lambda_{s}^{n}$.

\begin{definition}
For a given $\gamma_{t}\in\Lambda^{n}$ if we have
\[
v(\gamma_{t,s})=v(\gamma_{t})+a(s-t)+o(|s-t|)\ \text{as}\ s\rightarrow
t,\ s\geq t,\ \
\]
then we say that $v(\gamma_{t})$ is (horizontally) differentiable in $t$ at
$\gamma_{t}$ and denote $D_{t}v(\gamma_{t})=a$. $u$ is said to be horizontally
differentiable in $\Lambda^{n}$ if $D_{t}v(\gamma_{t})$ exists for each
$\gamma_{t}\in\Lambda^{n}$.
\end{definition}

\begin{definition}
Define $\mathbb{C}^{j,k}(\Lambda^{n})$ as the set of function $v:=(v(\gamma
_{t}))_{\gamma_{t}\in\Lambda^{n}}$ defined on $\Lambda^{n}$ which are $j$
times horizontally and $k$ times vertically differentiable in $\Lambda^{n}$
such that all these derivatives are $\Lambda^{n}$--continuous.
\end{definition}

The following It\^{o} formula was firstly obtained by Dupire \cite{Dupire.B}
and then by Cont and Fourni\'{e} \cite{Cont.R} for a more general formulation.

\begin{theorem}
[Functional It\^{o}'s formula]\label{w2} Let $(\Omega,\mathcal{F}%
,(\mathcal{F}_{t})_{t\in\lbrack0,T]},P)$ be a probability space, if $X$ is a
continuous semi-martingale and $u$ is in $\mathbb{C}^{1,2}(\Lambda^{n})$, then
for any $t\in\lbrack0,T)$,
\[
v(X_{t})-v(X_{0})=\int_{0}^{t}D_{s}v(X_{s})\,ds+\int_{0}^{t}D_{x}%
v(X_{s})\,dX(s)+\frac{1}{2}\int_{0}^{t}D_{xx}v(X_{s})\,d\langle X\rangle
(s),\quad\quad P-a.s..
\]

\end{theorem}

Note that the coefficients of (\ref{fbasd-1}) and (\ref{fbase-2}) are
permitted to be random. For a given $t\in\lbrack0,T]$ and $\gamma_{t}%
=(\gamma_{t}^{1},\gamma_{t}^{2})\in C^{d}\times C^{n}$, consider\ the
following functional fully coupled FBSDE:%
\begin{equation}%
\begin{array}
[c]{rl}%
W^{\gamma_{t}}(s)= & \gamma_{t}^{1}(t)+W(s)-W(t),\text{ }s\in\lbrack t,T],\\
X^{\gamma_{t}}(s)= & \gamma_{t}^{2}(t)+\int_{t}^{s}b(W_{r}^{\gamma_{t}}%
,X_{r}^{\gamma_{t}},Y^{\gamma_{t}}(r),Z^{\gamma_{t}}(r))dr+\int_{t}^{s}%
\sigma(W_{r}^{\gamma_{t}},X_{r}^{\gamma_{t}},Y^{\gamma_{t}}(r),Z^{\gamma_{t}%
}(r))dW(r),\text{ }s\in\lbrack t,T],
\end{array}
\label{path-1}%
\end{equation}%
\begin{equation}
Y^{\gamma_{t}}(s)=g(W_{T}^{\gamma_{t}},X^{\gamma_{t}}(T))-\int_{s}^{T}%
h(W_{r}^{\gamma_{t}},X_{r}^{\gamma_{t}},Y^{\gamma_{t}}(r),Z^{\gamma_{t}%
}(r))dr-\int_{s}^{T}Z^{\gamma_{t}}(r)dW(r),\quad s\in\lbrack t,T],
\label{path-2}%
\end{equation}
and%
\[%
\begin{array}
[c]{rl}%
W_{t}^{\gamma_{t}}(s)= & \gamma_{t}^{1}(s),\text{ \ \ }0\leq s\leq t,\\
X_{t}^{\gamma_{t}}(s)= & \gamma_{t}^{2}(s),\text{ \ \ }0\leq s\leq t.
\end{array}
\]
In order to apply Theorem \ref{theorem-1}, the above equations (\ref{path-1})
and (\ref{path-2}) can be rewritten as:%
\begin{equation}%
\begin{array}
[c]{rl}%
\tilde{X}^{\gamma_{t}}(s)= & \gamma_{t}(t)+\int_{t}^{s}\tilde{b}(\tilde{X}%
_{r}^{\gamma_{t}},Y^{\gamma_{t}}(r),Z^{\gamma_{t}}(r))dr+\int_{t}^{s}%
\tilde{\sigma}(\tilde{X}_{r}^{\gamma_{t}},Y^{\gamma_{t}}(r),Z^{\gamma_{t}%
}(r))dW(r),\\
\tilde{X}_{t}^{\gamma_{t}}= & \gamma_{t},\text{ \ \ }\tilde{X}_{s}^{\gamma
_{t}}=(W_{s}^{\gamma_{t}},X_{s}^{\gamma_{t}}),\text{ }t\leq s\leq T,
\end{array}
\label{path-11}%
\end{equation}

\begin{equation}
Y^{\gamma_{t}}(s)=g(W_{T}^{\gamma_{t}},X^{\gamma_{t}}(T))-\int_{s}^{T}%
h(\tilde{X}_{r}^{\gamma_{t}},Y^{\gamma_{t}}(r),Z^{\gamma_{t}}(r))dr-\int%
_{s}^{T}Z^{\gamma_{t}}(r)dW(r),\quad s\in\lbrack t,T], \label{path-21}%
\end{equation}
where $\tilde{b}$ and $\tilde{\sigma}$ are defined as%
\[%
\begin{array}
[c]{cl}%
\tilde{b}(\tilde{x}_{t}^{1},u^{1}(t))= & ((0,\cdots,0),b(\tilde{x}_{t}%
^{1},u^{1}(t)))^{T},\text{ \ \ }P-a.s.,\\
\tilde{\sigma}(\tilde{x}_{t}^{1},u^{1}(t))= & (I_{d\times d},\sigma(\tilde
{x}_{t}^{1},u^{1}(t)))^{T},\text{ \ \ }P-a.s.,
\end{array}
\]
$\forall(\tilde{x}_{T}^{1},u_{T}^{1}),(\tilde{x}_{T}^{2},u_{T}^{2})\in
C^{d+n}\times C^{m}\times C^{m\times d}$ and $(0,\cdots,0)\in\mathbb{R}^{d},$
$I_{d\times d}\ $is the $d$-dimensional identity matrix.

Now we consider the related P-PDE:
\begin{align}
&  D_{t}u(\gamma_{t})+\mathcal{L}u(\gamma_{t})-h(\gamma_{t},u(\gamma
_{t}),v(\gamma_{t}))=0,\label{path-3}\\
&  v(\gamma_{t})=D_{x}u(\gamma_{t})\tilde{\sigma}(\gamma_{t},u(\gamma
_{t}),v(\gamma_{t})),\nonumber\\
&  u(\gamma_{T})=g(\gamma_{T}^{1},\gamma^{2}(T)),\quad\gamma_{T}^{1}\in
\Lambda^{d},\text{ }\gamma^{2}(T)\in\mathbb{R}^{n},\nonumber
\end{align}
where
\[
\mathcal{L}u=(\mathcal{L}u_{1},\cdots,\mathcal{L}u_{n}),\quad\mathcal{L}%
=\frac{1}{2}\sum_{i,j=1}^{n+d}(\tilde{\sigma}\tilde{\sigma}^{T})_{i,j}%
(\gamma_{t},u,v)D_{x_{i}x_{j}}+\sum_{i=1}^{n+d}\tilde{b}_{i}(\gamma
_{t},u,v)D_{x_{i}}.
\]
\ \ 

\begin{theorem}
Suppose that Assumptions (\ref{ass-1}) and (\ref{ass-2}) hold. If $u^{1}%
\in\mathbb{C}^{1.2}(\Lambda^{d+n})$ is the solution of equation (\ref{path-3})
and $(u^{1},v^{1})$ is uniformly Lipschitz continuous and has linear growth,
then we have $u^{1}(\gamma_{t})=Y^{\gamma_{t}}(t)$, for each $\gamma
_{t}=(\gamma_{t}^{1},\gamma_{t}^{2})\in C^{d}\times C^{n},$ where
$(Y^{\gamma_{t}}(s),Z^{\gamma_{t}}(s))_{t\leq s\leq T}$ is the unique solution
of equations (\ref{path-11}) and (\ref{path-21}). Consequently,\ the P-PDE
(\ref{path-3})\ has a\ unique $\mathbb{C}^{1.2}(\Lambda^{d+n})$-solution.
\end{theorem}

\textbf{Proof}. Note that $(u^{1},v^{1})$ is uniformly Lipschitz continuous
and has linear growth. Then, under Assumptions (\ref{ass-1}) and
(\ref{ass-2}), the equation
\begin{align*}
&  d\tilde{X}^{\gamma_{t}}(s)=\tilde{b}(\tilde{X}_{s}^{\gamma_{t}}%
,u^{1}(\tilde{X}{_{s}^{\gamma_{t}}}),v^{1}(\tilde{X}{_{s}^{\gamma_{t}}%
}))ds+\tilde{\sigma}(\tilde{X}_{s}^{\gamma_{t}},u^{1}(\tilde{X}{_{s}%
^{\gamma_{t}}}),v^{1}(\tilde{X}{_{s}^{\gamma_{t}}}))dW(s),\\
&  \tilde{X}_{t}^{\gamma_{t}}=\gamma_{t},\quad s\in\lbrack t,T],
\end{align*}
has a unique solution. It is easy to see that we can denote $\tilde{X}%
_{s}^{\gamma_{t}}$ by $(W_{s}^{\gamma_{t}},X_{s}^{\gamma_{t}})$.

Set $(Y^{1}(s),Z^{1}(s))=(u^{1}(\tilde{X}{_{s}^{\gamma_{t}}}),v^{1}(\tilde
{X}{_{s}^{\gamma_{t}}})),\quad t\leq s\leq T$. Applying\ the\ functional
It\^{o} formula to $u^{1}({\tilde{X}_{s}^{\gamma_{t}}})$, we have
\begin{align*}
&  dY^{1}(s)=h(\tilde{X}_{s}^{\gamma_{t}},Y^{1}(s),Z^{1}(s))dr+Z^{1}%
(s)dW(s),\\
&  Y^{1}(T)=g(W_{T}^{\gamma_{t}},X^{\gamma_{t}}(T)),\quad\ \ s\in\lbrack t,T].
\end{align*}
Thus, by Theorem \ref{theorem-1}, we have $(Y^{1}(s),Z^{1}(s))=(Y^{\gamma_{t}%
}(s),Z^{\gamma_{t}}(s))$. In particular, $u^{1}(\gamma_{t})=Y^{\gamma_{t}}%
(t)$. This completes the proof. $\ \ \Box$
%%%%%%%%%%%%%%%%%%%%%%%²Î¿¼ÎÄÏ×
\renewcommand{\refname}{\large References}\setlength{\itemsep}{-0.0em}

\bigskip

\textbf{Acknowledgements} The authors would like to thank Prof. Shige Peng for
some useful conversations.

\end{document}